\documentclass[11pt]{article}
\input epsf.tex
\usepackage{latexsym}
\usepackage{amsmath,amsthm,amsfonts,amssymb}
\usepackage{epsfig}

\newtheorem{pro}{Proposition}[section]
  \newtheorem{thm}[pro]{Theorem}
  \newtheorem{lem}[pro]{Lemma}

  \newtheorem{cor}[pro]{Corollary}

  \def\C{{\cal C}}
  \def\H{{\cal H}}
  
  \def\P{{\cal P}}
  \def\O{{\cal O}}
  \def\M{{\cal M}}
  \def\F{{\cal F}}
  \def\T{{\cal T}}

  \def\chix{{\raise.5ex\hbox{$\chi$}}}

\begin{document}
\title{On the existence of completely saturated packings and completely
reduced coverings}
\author{Lewis Bowen}
\maketitle
\begin{abstract}
We prove the following conjecture of G. Fejes Toth, G. Kuperberg, and
W. Kuperberg: every body $K$ in either $n$-dimensional Euclidean or
$n$-dimensional hyperbolic space admits a completely saturated packing
and
a completely reduced covering. Also we prove the following
counterintuitive result: for every $\epsilon > 0$, there is a body $K$
in
hyperbolic $n$-space which admits a completely saturated packing with
density less than $\epsilon$ but which also admits a
tiling. \end{abstract}

\section{Introduction}

    Two of the most basic problems in the theory of packings and
coverings
is to find the most efficient packing and covering by replicas of a
given
set $K$ in either $n$-dimensional Euclidean space or $n$-dimensional
hyperbolic space ([FeK]). The term `most efficient' is intentionally
ambiguous as different circumstances have led to different
interpretations. Generally, the density of a packing of Euclidean space
is
defined as the limit of the relative fraction of volume occupied by the
bodies in the packing in a ball of radius $r$ centered at some point $p$
as $r$ tends to infinity. It can be shown that if this limit exists for
some point $p$, it exists for every point $p$ and has the same value. A
densest packing is then one which maximizes this density. It is easy to 
prove the existence of a densest packing in Euclidean space. The density of 
a
covering can be defined similarly and a thinnest covering is one that
minimizes the density.

One disadvantage of this definition is that is does not depend on local
structure. For example, if a finite number of bodies in a densest
packing
are removed, it remains a densest packing. Another disadvantage is that
it
is not clear that this definition carries over well to hyperbolic
space. For instance, the above limit may exist for every point in
hyperbolic space yet take on different values for different points. For
these and other reasons, the authors of [FKK] introduced the concept of
a
completely saturated packing. For such a packing, it is not possible to
replace a finite number of bodies of the packing with a greater number
of
bodies and still remain a packing. Similarly a completely reduced
covering
is one in which it is not possible to replace a finite number of bodies
of
the covering with a smaller number of bodies and still remain a
covering.

In [FKK], it was proven that any convex body of Euclidean space admits a
completely saturated packing (and more generally any body with the
strict
nested similarity property) (see also [Kup]). The authors conjectured
that
completely saturated packings and completely reduced coverings exist for
any body $K$ in either $n$-dimensional Euclidean or hyperbolic space. In
section 4, we prove this conjecture. The proof requires few properties of 
either Euclidean or hyperbolic space and easily extends to a more general 
setting (see Remark after Theorem 2.2).

    Our techniques are based on those developed in [BoR]. In that paper,
a
new notion of efficiency was introduced which intuitively
applies
to a class of measures on a space of packings rather than to individual
packings. Such measures are often easier to analyze but at the same time
they carry information about the `average' properties of some packings.
In
particular, Birkhoff's ergodic theorem (in the Euclidean case) and
recent
results due to Nevo and Stein (in the hyperbolic case) show that this
new
notion (to be explained in section 3) coincides (in a measure-theoretic
way) to the usual notion of density given by limits of volumes in
expanding balls. This will enable us to construct a peculiar family of
bodies $K$ that admit tilings of hyperbolic space but such that they
also
admit completely saturated packings with arbitrarily low density.

{\bf Acknowledgments}:  I would like to thank my advisor Charles Radin
for
many useful conversations and suggestions.

\section{Notation and Results}

Let $S$ denote either $n$-dimensional hyperbolic or $n$-dimensional 
Euclidean space with a
distinguished point $\O$ which we call the origin. Let $G$ be the group of 
orientation-preserving isometries of $S$.

  A \textbf{body} $K$ is a nonempty connected compact subset of $S$ which
equals the closure of its interior. Let $\Sigma_K = \{g \in G|gK=K\}$.

A \textbf{packing} by a body $K$ is a collection $\P=\{g_1K,g_2K,...\}$
of
congruent copies of $K$ that pairwise have disjoint interiors. Note that
if $g_1K \in \P$ for some $g_1 \in G$, then $g_1k \in \P$ for all $k \in
\Sigma_K$ since $kK = K$ implies $g_1kK = g_1K$. Conversely, if $g_1K =
g_2K$ then $g_2^{-1}g_1 \in \Sigma_K$. So we let $\hat \P$ be the subset
of the space $G/\Sigma_K$ of left cosets defined by $g\Sigma_K \in \hat
\P$ if any only $gK \in \P$. We note that $\hat \P$ uniquely determines
$\P$ (given $K$) so if $\F \subset G/\Sigma_K$, we will say that $\F$ 
determines the
packing $\P$ if $\hat \P = \F$. Later, we will use this correspondence to
define a metric on the space of packings.

A \textbf{covering} by a body $K$ is a collection $\C=\{g_1K,g_2K,...\}$
of congruent copies of $K$ such that for every point $p \in S$, $p \in
g_jK$ for some $g_jK \in \C$. As above, we let $\hat \C$ be the subset
of
the space $G/\Sigma_K$ of left cosets defined by $g\Sigma_K \in \hat \C$
if any only $gK \in \C$. A \textbf{tiling} by $K$ is a packing that is
also a covering.

  A \textbf{completely saturated} packing $\P$ is such that for no packing
$\F_1 \subset \P$ does there exist another finite packing $\F_2$ with more 
bodies that $\F_1$ such that $(\P -
\F_1) \cup \F_2$ is a packing.

A \textbf{completely reduced} covering $\C$ is such that for no finite
subset $\F_1 \subset \hat \C$ does there exist another finite set $\F_2
\subset G/\Sigma_K$ with fewer elements that $\F_1$ such that $(\hat \C
-
\F_1) \cup \F_2$ determines a covering. Note that this is the equivalent to 
the definition given earlier.

Our main result is the following.

  \begin{thm}
  For any body $K$ in $S$, there exists a completely saturated packing by
  $K$ and a completely reduced covering by $K$.
  \end{thm}

\noindent {\bf Remark}: The theorem extends with almost no changes to the 
case in which $S = G/H$ where $G$ is a unimodular Lie group, $H$ is a 
compact subgroup of $H$ and $G$ satisfies the following property. For every 
compact subset $C$ of $G$ there exists a discrete cocompact subgroup $G' < 
G$ such that there exists a fundamental domain $F$ for $G'$ and $C \subset 
F$. In this setting, `congruent copies of a body $K \subset S$' means 
subsets of the form $gK$ where $g \in G$ and $G$ acts on subsets of $S$ in 
the natural way.\\
\\

We will only prove the existence of completely saturated packings for
$K$
as the other assertion can be proven similarly.

For any point $p \in S$, we let $B_n(p)$ denote the closed ball of radius 
$n$
centered at $p$. Also we let $\lambda_S$ denote the usual measure on
$S$. Finally for any packing $\P$ we let $c(\P) = \cup_{gK \in \P} \,
gK$
be the closed subset associated to $\P$.

We obtain the following
counterintuitive result:

\begin{thm} Let $S$ be $n$-dimensional hyperbolic space for some
$n>1$. Let $\epsilon > 0$. Then there is a body $K$ in $S$ that admits a
tiling $\T$ and a completely saturated packing $\P$ with the following
properties. For every point $p \in S$, the limit
\begin{equation}
\lim_{r \to \infty} \frac{\lambda_S(c(\P) \cap
B_r(p))}{\lambda_S(B_r(p))}
\end{equation}
exists and is independent of $p$. Furthermore this common limit is less
than $\epsilon$.
\end{thm}

The paper is organized as follows. In section 3, we introduce the main
tools of our study. In section 4, we prove Theorem 2.1. In section 5,
we
prove Theorem 2.2.

  \section{Measures on the Space of Packings}

Let $\Sigma_\O$ be the stabilizer $\{g \in G| \, g\O = \O\}$ of the origin 
$\O$. It is isomorphic to $SO(n-1)$. Let $\pi: G \to S$ be the
canonical projection identifying $G/\Sigma_\O$ with $S$ such that $\pi(g) = 
g\O$ for every $g \in G$. Let $\lambda_G$ be the normalization of Haar 
measure on
$G$ such that $\lambda_S$, the usual measure on $S$, has the following 
relationship with $\lambda_G$ (see [Rat] for the hyperbolic case). For any 
Borel set $E \subset S$,
$\lambda_S(E) = \lambda_G(\pi^{-1}(E))$. We note the general fact (see
[Lan]) that the Haar measure on any locally compact group $G_0$ is
determined uniquely up to scalar multiplication by the property that it
is
a nonzero Borel measure that is invariant under left multiplication by
every element $g \in G_0$ (i.e. $\lambda_{G_0}(E) = \lambda_{G_0}(gE)$
for
every $g\in G$ and Borel set $E \subset G$ where $\lambda_{G_0}$ is a
Haar
measure of $G_0$).

Let $d_S$ be the usual distance function on $S$. Let $d_G$ be any 
left-invariant distance function on $G$ satisfying $d_S(p_1,p_2) = \inf 
\{d_G(\tilde p_1, \tilde p_2)  | \, \tilde p_1 \in
\pi^{-1}(p_1)$ and $\tilde p_2 \in \pi^{-1}(p_2)\}$. The left-invariant of 
$d_G$ means that for any $g_0, g_1, g_2 \in G$, $d_G(g_1,g_2) = d_G(g_0 
g_1,g_0 g_2)$. For any body $K$ in
$S$, we define a distance function $d^K$ on $G/\Sigma_K$ by
$d^K(g_1\Sigma_K,g_2\Sigma_K) = \inf \{d_G(g_1k_1,g_2k_2)  | \,
k_1,k_2\in\Sigma_K\}$. In other words, $d^K$ is the pushforward of
$d_G$.

Let $B_n, B^K_n,$ and $\tilde B_n$ denote the closed balls of radius $n$
centered at $\O$ in $S$, $\Sigma_K$ in $G/\Sigma_K$ and the identity in $G$
respectively. If $X$ is a subset of $G/\Sigma_K$ and $\epsilon > 0$
we let $N_\epsilon(X)$
  denote the open $\epsilon$-neighborhood of $X$ in $G/\Sigma_K$.

  Let $C_K$ denote the set of packings of $S$ by $K$. First, for every $n
 >
0$, define a pseudometric $d_n$ on $C_K$ by
  \begin{equation}
  d_n(\P_1,\P_2) = \inf \{\epsilon > 0 | \, B^K_n \cap \hat \P_1 \subset
  N_\epsilon(\hat \P_2) \textnormal{ and } B^K_n \cap \hat \P_2 \subset
N_\epsilon(\hat \P_1)
  \}.
  \end{equation}
  \noindent Next we define a metric $d_K$ on $C_K$ by:
  \begin{equation}
  d_K(\P_1,\P_2) = \sup_{n > 0}
\min\{\frac{d_n(\P_1,\P_2)}{n},\frac{1}{n}\}.
  \end{equation}

  \noindent It can be checked that $C_K$ is a compact metric space on
which
$G$ acts
  jointly continuously in the obvious way, i.e. for $\P \in C_K$ and $h
\in
G$, $h\P = \{hgK | \, gK \in \P\}$ (cf [RaW]). We note that in Euclidean
space, if $h$ is translation by $x$, then $h\P$ is commonly denoted by
$\P
+ x$. Although the metric $d_K$ will be convenient in what is to follow,
we are really only concerned with the topology that it induces. The idea
is that two packings are close in $C_K$ if in a large ball about the
origin, the bodies of one packing can by small rigid motions be made to 
coincide
with the bodies of the other packing.

Let $M=M(K)$ be the set of Borel
  probability measures on $C_K$. We let $M_I = M_I(K)$ be the
  subset of $M$ consisting of all those measures that are $G$-invariant.
In
other words, $\mu \in M$ is in $M_I$ if and only if $\mu(gE) = \mu(E)$
for
every $g \in G$ and Borel set $E \subset C_K$.

Although we will not use it, we show here how to construct a natural
class
of measures in $M_I$. If $\P$ is any packing such that the symmetry
group
$\Sigma_\P = \{g \in G| \, g\P = \P\}$ of $\P$ is such that $G/\Sigma_K$
is compact then there is a canonical measure $\mu_\P \in M_I$ associated
to $\P$. There is a homeomorphism $q$ from the orbit $O(\P) =\{g\P | \,
g
\in G\}$ to $G/\Sigma_\P$ defined by $q(g\P) = g\Sigma_\P$. This
homeomorphism commutes with the action by $G$ where $G$ acts on
$G/\Sigma_\P$ by left multiplication. There is a natural measure on
$G/\Sigma_\P$ induced from $\lambda_G$ since the quotient map from $G$
to
$G/\Sigma_\P$ is a covering map in which the covering transformations
are
all measure-preserving (in fact they are isometries of $G$). This
measure
is invariant under the action of $G$ since $\lambda_G$ is invariant. Now
a
measure $\hat \mu_\P$ on $O(\P)$ can be induced (via $q$) from the
measure
on $G/\Sigma_\P$. This measure can then be extended to a measure
$\mu_\P$
on all of $C_K$ by $\mu_\P(E) = \hat \mu_\P(E \cap O(\P))$. Thus for
every
packing whose symmetry group is cocompact, there is a natural measure in
$M_I$ associated to that packing.

  For any packing $\P$ by a body $K$, we let $c(\P) = \cup_{gK \in \P} \,
gK \subset S$ be the closed subset associated to $\P$. The density of a
measure $\mu \in M_I$, $D(\mu)$ is defined by
  \begin{equation}
  D(\mu) := \mu(\{\P \in C_K| \, \O \in c(\P)\}).
  \end{equation}

\noindent In other words, the density of $\mu$ is the probability (with
respect to $\mu$) that the origin lies in a body of a random packing. We
define the (measure-theoretic) \textbf{optimal density} of a body $K$
  to be $D(K) := \sup_{\mu \in M_I} D(\mu)$. A measure $\mu \in M_I$ is
  said
  to be \textbf{optimally dense} (for $K$) if $D(\mu) = D(K)$ and $\mu$
is
ergodic. By ergodic we mean that for any $E \subset C_K$ such that $gE =
E$ for all $g \in G$, $\mu(E) \in \{0,1\}$.

We will sketch a proof that when $S$ is Euclidean, the optimal density
of
$K$ is equal to the usual notion. To be precise, for a packing $\P$ and
any point $p \in S$, let $D(\P,p)$ be defined by
\begin{equation}
\delta(\P,p) = \lim_{r\to\infty} \frac{\lambda_S(c(\P)\cap
B_r(p)}{\lambda_S(B_r)}
\end{equation}
\noindent whenever this limit exists. Let $\delta(K) = \sup_{\P \in C_K}
\delta(\P,p)$. Note that $\delta(K)$ does not depend on the point $p$ chosen
since $\delta(\P,p) = \delta(g\P,\O)$ where $g \in G$ such that $gp=\O$. We 
will
show that $\delta(K) = D(K)$.

Groemer [Gro] proved that, when $S$ is Euclidean, for any body $K
\subset
S$ there exists a packing $\P$ by $K$ such that $\delta(\P,p) = \delta(K)$
for
all points $p \in S$ and such that this limit converges uniformly in
$p$. By standard ergodic theory, there exists a measure $\mu \in M_I$
whose support is contained in the closure of the orbit $O(\P) = \{g\P |
\,
g \in G\}$. Since convergence is uniform in $p$, $\delta(\P',\O) = 
\delta(\P,\O)$
for all $\P'$ in the closure of $O(\P)$. By Birkhoff's ergodic theorem,
$\delta(\P',\O) = D(\mu)$ for $\mu$-almost every packing $\P'$. Thus $D(\mu)
=
\delta(K)$ and so $D(K) \ge \delta(K)$.

For the other way, let $\mu$ be an optimally dense measure (the
existence
of which is proven below). Then Birkhoff's ergodic theorem implies that
for $\mu$-almost every packing $\P$, $\delta(\P,\O) = D(\mu)$. Thus
$\delta(K) \ge D(K)$. So the two notions are equal. In the hyperbolic
case, we can only show that $\delta(K) \ge D(K)$ (see Corollary 5.2). In
fact,
for any given $\epsilon$ we will construct a body $K$ such that there
exists a tiling by $K$ (and hence $\delta(K) = 1$) but $D(K) <
\epsilon$. This immediately implies that the subspace of $C_K$
consisting
of all tilings by $K$ has measure zero with respect to any measure in
$M_I$. So, by restricting our attention to invariant measures, many
`pathological' packings (or tilings) are ignored. It may seem odd to some 
readers that we call $D(K)$ the optimal density in spite of the fact 
$\delta(K)$ may be greater than $D(K)$. So we emphasize that $D(K)$ is just 
the {\it measure theoretic} optimal density.

\noindent The main theorem of this paper is:

  \begin{thm}
  Let $U \subset C_K$ be the set of all packings by $K$ that are not
  completely saturated. Then for any optimally dense measure $\mu$ for
  $K$,
  $\mu(U) = 0$.
  \end{thm}

\noindent Theorem 2.1 follows immediately from this one and the
following
existence theorem.

\begin{thm}
For any body $K \subset S$, there exists an optimally dense measure
$\mu$
for $K$.
\end{thm}
\noindent A proof of this theorem appears in [BoR, Theorem 3.2]. We include 
the same
proof here for completeness.

\begin{proof}

We let $M_I$ have the weak* topology. Equivalently $\mu_n \to \mu \in
M_I$
if and only if for every continuous function $f: C_K \to {\mathbb R}$,
$\int_{C_K} f \, d\mu_n \to \int_{C_K} f \, d\mu$. By standard functional
analysis, since $C_K$ is a compact metric space, $M_I$ is compact.

Let $Z =\{\P \in C_K | \, \O \in c(\P)\}$. Let $\chi_Z$ be the
characteristic function of $Z$. Because $\chi_Z$ is upper semicontinuous
there exists a
      decreasing sequence $f_j$ of continuous real valued functions on
   $C_K$
      which converge pointwise to $\chi_Z$. Choose a sequence $\mu_k\in
      \M_I$
      such that $D(\mu_k)=\int_{C_K} F_\O \, d\mu_k\to D(K)$ as $k\to
      \infty$, and, using the compactness of $\M_I$, assume without
  loss
     of
      generality that $\mu_k$ converges to some $\mu_\infty\in \M_I$.
    Then
      $\int_{C_K} f_j\, d\mu_k\to \int_{C_K} f_j\, d\mu_\infty$ as $j\to
      \infty$, and $\int_{C_K} f_j\, d\mu_\infty\searrow D(\mu_\infty)$
as
      $k\to \infty$. Since $\int_{C_K} f_j\, d\mu_k\ge D(\mu_k)$ and
     $D(\mu_k) \to D(K)$ as $k\to \infty$, $D(\mu_\infty)\ge D(K)$. From
   the
      Krein-Milman
      theorem there exists an ergodic measure $\tilde \mu \in \M_I$ for
which $D(\tilde
      \mu)=\int_{C_K} \chi_Z \, d\tilde\mu \ge D(\mu_\infty)$, and thus
      $D(\tilde\mu)\ge D(K)$. But then from the definition of $D(K)$,
      $D(\tilde\mu)= D(K)$.
\end{proof}

  \section{A Family of saturating maps}

     Before proving the theorem, we outline the case in which $S$ is the
  Euclidean plane. Let $T_j$ be the usual square tiling by
  squares of side length $j$. We assume that the origin is in the center
  of a tile of $T_j$
  for all $j$. Let $G_j$ be the group of translations which fix
  $T_j$. Let $U_j$ be the set of packings which are not saturated in the
square of $T_j$ that contains the origin. To be precise, $\P \in U_j$ if
and only if there is a finite packing $\F_1 \subset \P$ and a finite
packing $\F_2$ such that $\F_2$ has more elements than $\F_1$, for all
$gK
\in \F_1 \cup \F_2$, $gK$ is contained in the square of $T_j$ that
contains the origin and $(\P - \F_1) \cup \F_2$ is a packing.

For each $j$, there exists a Borel map $\Phi_j: C_K \to C_K$
  satisfying:
\begin{enumerate}
\item $gK \in \P$ and $gK$ intersects an edge of $T_j$ if and only if
  $gK \in \Phi_j(\P)$ and $gK$ intersects an edge of $T_j$;
\item For any given square in $T_j$, the number of bodies of
$\Phi_j(\P)$
that are contained in that square is as large as possible given the
above
constraint;
\item $\Phi_j$ commutes with $G_j$.
\end{enumerate}
  \noindent The idea is that $\Phi_j$ replaces $\P \in C_K$ with a
$\Phi_j(\P)$ that
  is saturated with respect to the squares of $T_j$. Because $\Phi_j$
commutes with $G_j$, it induces a map $\Phi_{j*}:M \to M$ such that if
$\mu \in M_I$ then $\Phi_{j*}(\mu)$ is invariant under $G_j$. Since
$G_j$
is
  cocompact and discrete we can average $\Phi_{j*}(\mu)$ over a
  fundamental
  domain for $G_j$ in $G$ to obtain a measure that is invariant under all
  of
  $G$.
This new measure will have a density greater than $\mu$ unless
  $\mu(U_j) = 0$. Since $\cup_j U_j = U$ is the set of all non-completely
saturated packings, this shows that an optimally dense measure must
satisfy $\mu(U) = 0$ which proves the theorem.

  By a fundamental domain $F$ of a subgroup $H < G$, we shall mean a
connected set in $S$ equal to the closure of its interior such that $\{hF| 
\, h
\in H\}$ is a tiling by $F$.

For the general case, let $\{G_j\}_{j=0}^\infty$ be a sequence of
  discrete cocompact subgroups of $G$ such that there
  exist
  fundamental
  domains $F_j$ (in $S$) for $G_j$ with $B_j \subset F_j$. We will also
assume that $F_j$ and $G_j$ have been chosen so that for all $g \in
\pi^{-1}(F_j)$, $g^{-1} \in \pi^{-1}(F_j)$. For the Euclidean case, we 
could, for example, let $G_j$ be the cubic lattice generated by the 
translations $\tau_{i,j}:=(x_1,...,x_n) \to (x_1,...,x_i+j,x_{i+1},...,x_n)$ 
for $1\le i \le n$. Then we could choose $F_j$ to be the cube of side length 
$j$ whose center is the origin and whose faces are parallel to the 
coordinate planes. For the hyperbolic case, we refer to [FKK, Theorem 4.1] 
for the existence of $\{G_j\}_{j=0}^\infty$.

  For $\P \in C_K$ and $F \subset S$, let $\P*F$ be the packing
  consisting of all elements $gK \in \P$ such that the interior of $gK$
intersects $F$ nontrivially. Also let
  $\partial F$ denote the boundary of $F$, i.e. the intersection of $F$
with the closure of its complement. Let $X_j = \{x \in C_K | x* \partial
F_j = x\}$ with the subspace topology. For any $x \in X_j$, a filling
$f$
for $x$ is a packing by $K$ such that $f * \partial F_j = x$, $f * F_j =
f$ and the number of elements of $f$ is as large as possible given these
constraints.

  \begin{lem}
  For each $j \ge 0$, there exists a Borel map $\phi_j : X_j \to C_K$
  such that for each $x \in X_j$, $\phi_j(x)$ is a filling for $x$.
  \end{lem}

  \begin{proof}

Let $j \ge 0$ be fixed and let $X = X_j$. Given $x \in X$, $f$ a filling
for $x$ and $m > 0$, let $Y(f,m) = \{ x' \in X| $ there exists a filling
$f'$ for $x'$ with $D_K(f',f) < 2^{-m}\}$. From the compactness of $F_j$
it follows that given $m > 0$, there
  exists
  a finite number of packings $x_1,...,x_p \in X$ and fillings $f_k$ for
  $x_k$ such that $X = \bigcup_{k=1}^p Y(f_k,m)$.

  Thus for a given integer $m > 0$, there exists a finite partition
  $\{A_{m,k}\}_{k=1}^{r_m}$ of $X$ such that each $A_{m,k}$ is Borel and
contained in
  some $Y(f,m)$. We will assume that $\{A_{m+1,k}\}_{k=1}^{r_{m+1}}$
refines
  $\{A_{m,k'}\}_{k'=1}^{r_m}$ (i.e. for every $A_{m+1, k}$, there exists
an
$A_{m,k'}$ such that $A_{m+1,k} \subset A_{m,k'}$). For each $m > 0$ and
$1 \le k \le r_m $, choose $a(m,k) \in A_{m,k}$ and $f(m,k)$ a filling
for
$a(m,k)$ so that the following are satisfied:
\begin{enumerate}
\item $\{a(m,k)\}_{k=1}^{r_m} \subset \{a(m+1,k')\}_{k'=1}^{r_{m+1}}$
for
all $m > 0$.
\item For every $m > 0$, $1 \le k \le r_m$, $a \in A_{m,k}$ there is a
filling $f$ for $a$ such that $d_K(f,f(m,k)) < 2^{-m+1}$. This is
possible
since for each $A_{m,k}$ there exists $x \in X$ and a filling $f$ for
$x$
with $A_{m,k} \subset Y(f,m)$.
\item If for some $m$ and $k, k'$, $A_{m+1,k'} \subset A_{m,k}$, then
$d_K(f(m+1,k'),f(m,k)) < 2^{-m+1}$.
\end{enumerate}

  We define a map $\alpha_m: X \to C_K$ by $\alpha_m(x) = f(m,k)$ if $x
\in
  A_{m,k}$. Each $\alpha_m$ is Borel (since each $A_{m,k}$ is Borel). We
claim that $\{\alpha_m\}_{m=1}^\infty$ converges pointwise as to a map
$\phi$. Let $x \in X$. Then by the third condition above, if $m < m'$
then
\begin{equation}
D_K(\alpha_m(x), \alpha_{m'}(x)) \le \Sigma_{k=m}^{m'-1} \, 2^{-k+1} \le
\Sigma_{k=m}^\infty \, 2^{-k+1} \to 0
\end{equation}
\noindent as $m \to \infty$. So the sequence
$\{\alpha_m(x)\}_{m=1}^\infty$ is Cauchy and therefore converges in
$C_K$
to an element $\phi(x)$. Since all the maps $\alpha_m$ are Borel, $\phi$
must be Borel, too.

We claim that $\phi(x)$ is a filling for $x$ for each $x \in X$. So let
$x
\in X$ and for each $m > 0$ choose a filling $f_m$ for $x$ such that
$d_K(f_m, \alpha_m(x)) < 2^{-m+1}$. By the definition of $A_{m,k}$ and
the
second condition of $f(m,k)$ listed above, such an $f_m$ exists for all
$m
 >0$. By the triangle inequality, $d_K(f_m, \phi(x)) \le 2^{-m+1} +
\Sigma_{k=m}^\infty 2^{-k+1}$ which goes to zero as $m$ goes to
infinity. Hence the sequence $\{f_m\}_{m>0}$ converges to $\phi(x)$.
Since
$f_m$ is a filling for $x$ for each $m$, it follows that $\phi(x)$ is
also
a filling for $x$. Now we are done.

  \end{proof}
We choose maps $\Phi_j$ satisfying the conclusion of the above
lemma. Define $\Phi_j: C_K \to C_K$ by the following properties:
\begin{enumerate}
\item for any $\P \in C_K$, $\Phi_j(\P) * F_j = \phi_j(\P \cap \partial
F_j) * F_j$;
\item for any $g \in G_j$ and $\P \in C_K$, $\Phi_j (g\P) = g
\Phi_j(\P)$.
\end{enumerate}
  It should be clear from this definition that for each $j$, $\Phi_j$ is
Borel. Let $\mu \in M_I$ be given. Let $\mu'_j = \Phi_{j*}(\mu)$, i.e.
for
every Borel set $E \subset C_K$, $\mu'_j(E) = \mu(\Phi_j^{-1}(E))$. By
the
above,
  $\mu'_j$ is a Borel probability measure that is invariant under
$G_j$. Let $\mu_j \in M_I$ be defined by

  \begin{equation}
  \mu_j(E) = \frac{1}{\lambda_S(F_j)} \int_{\pi^{-1}(F_j)}
\mu'_j(g^{-1}E) \, d\lambda_G(g)
  \end{equation}

  \noindent for any Borel $E \subset C_K$. It is easy to show that
$\mu_j$
is $G$-invariant and therefore
  really is in $M_I$. The next two lemmas will provide tools for
calculating $D(\mu)$ and $D(\mu_j)$ with respect to the relative density
within $F_j$.

  \begin{lem}
  \begin{equation}
  D(\mu_j) = \int_{C_K} \frac{\lambda_S(c(\P) \cap F_j)}{\lambda_S(F_j)}
\,
d\mu'_j(\P)
  \end{equation}
  \end{lem}

  \begin{proof}

  Let $Z= \{\P \in C_K| \O \in c(\P)\}$. Let $\chi_Z$ denote the
characteristic function of $Z$. By definition of density and of $\mu_j$,

  \begin{eqnarray}
  D(\mu_j) &=& \mu_j(Z)\\
           &=& \frac{1}{ \lambda_S(F_j)} \int_{\pi^{-1}(F_j)}
\mu'_j(g^{-1}Z)\,
  d\lambda_G(g)\\
           &=& \frac{1}{\lambda_S(F_j)} \int_{\pi^{-1}(F_j)} \int_{C_K}
  \chi_{g^{-1}Z}(\P) \, d\mu'_j(\P) \, d\lambda_G(g)\\
           &=& \int_{C_K} \frac{1}{\lambda_S(F_j)} \int_{\pi^{-1}(F_j)}
  \chi_{g^{-1}Z}(\P) \, d\lambda_G(g) \, d\mu'_j(\P)
  \end{eqnarray}

\noindent Hence we will be done once we show that for any $\P \in C_K$,
\begin{equation}
\int_{\pi^{-1}(F_j)}
  \chi_{g^{-1}Z}(\P) \, d\lambda_G(g) = \lambda_S(c(\P) \cap F_j).
\end{equation}

\noindent So,

\begin{eqnarray}
\int_{\pi^{-1}(F_j)}
  \chi_{g^{-1}Z}(\P) \, d\lambda_G(g) &=& \lambda_G(\{g \in
\pi^{-1}(F_j)|
\, \P \in g^{-1}Z\})\\
                               &=& \lambda_G(\{g \in \pi^{-1}(F_j)| \,
g\P
\in Z\})\\
                               &=& \lambda_G(\{g \in \pi^{-1}(F_j)| \, \O
\in c(g\P)\})\\
                               &=& \lambda_G(\{g \in \pi^{-1}(F_j)| \,
g^{-1}\O \in c(\P)\})\\
                               &=& \lambda_G(\{g \in \pi^{-1}(F_j)| \,
g\O
\in c(\P)\})\\
                               &=& \lambda_G(\pi^{-1}(c(\P) \cap F_j))\\
                               &=& \lambda_S(c(\P) \cap F_j)
\end{eqnarray}

Equation (18) holds because $\lambda_G$ is inversion-invariant (see [Rat] 
for the hyperbolic case) and because $\pi^{-1}(F_j) =
\pi^{-1}(F_j)^{-1}$. The last equation holds since for any Borel set $E 
\subset S$, $\lambda_S(E) =
\lambda_G(\pi^{-1}(E))$.

  \end{proof}

  \begin{lem}
  \begin{equation}
  D(\mu) = \int_{C_K} \frac{\lambda_S(c(\P) \cap F_j)}{\lambda_S(F_j)} \,
d\mu
  \end{equation}
  \end{lem}
  \begin{proof} By equation (13) in the previous lemma,
  \begin{eqnarray}
  \int_{C_K} \frac{\lambda_S(c(\P) \cap F_j)}{\lambda_S(F_j)} \,
d\mu(\P) &=& \int_{C_K}
  \int_{\pi^{-1}(F_j)} \frac{\chi_Z(g^{-1}\P)}{\lambda_S(F_j)} \,
d\lambda_G(g) \, d\mu(\P)\\
     &=&  \int_{\pi^{-1}(F_j)} \int_{C_K}
\frac{\chi_Z(g^{-1}\P)}{\lambda_S(F_j)}
  \, d\mu(\P) \, d\lambda_G(g)\\
     &=& \int_{\pi^{-1}(F_j)} \int_{C_K}
\frac{\chi_Z(\P)}{\lambda_S(F_j)}
\,
  d\mu(\P)
  \, d\lambda_G(g)\\
     &=& \int_{\pi^{-1}(F_j)} \frac{D(\mu)}{\lambda_S(F_j)}
d\lambda_G(g)\\
     &=& D(\mu)
  \end{eqnarray}

  \noindent The third inequality holds because $\mu$ is $G$-invariant.
  \end{proof}
\begin{proof}{\it (of Theorem 3.1)} Let $U_j$ be the set of all packings
in $C_K$ that are unsaturated
  relative to $F_j$, i.e. for any $\P \in U_j$,
  there is a packing $\P'$ such
  that $\P * \partial F_j = \P'* \partial F_j$ but $\P' * F_j$ has more
elements than $\P * F_j$.

  \begin{eqnarray}
  D(\mu_j) &=& \int_{C_K} \frac{\lambda_S(c(\P) \cap
F_j)}{\lambda_S(F_j)}\, d\mu'_j(\P)\\
           &=& \int_{C_K} \frac{\lambda_S(c(\Phi_j(\P)) \cap
F_j)}{\lambda_S(F_j)} \,
  d\mu(\P)\\
           &=& \int_{U_j} \frac{\lambda_S(c(\Phi_j(\P)) \cap
  F_j)}{\lambda_S(F_j)}\, d\mu(\P)\\
           &+& \int_{C_K - U_j} \frac{\lambda_S(c(\Phi_j(\P)) \cap
  F_j)}{\lambda_S(F_j)}\, d\mu(\P)\\
          &\ge& \int_{U_j} \frac{\lambda_S(c(\P) \cap F_j) +
\lambda_S(K)}{\lambda_S(F_j)}\,
  d\mu(\P)\\
          &+& \int_{C_K - U_j} \frac{\lambda_S(c(\Phi_j(\P)) \cap
F_j)}{\lambda_S(F_j)}\,
  d\mu(\P)\\
           &=& \int_{C_K} \frac{\lambda_S(c(\P) \cap
F_j)}{\lambda_S(F_j)}
\, d\mu(\P)+
  \mu(U_j)\frac{\lambda_S(K)}{\lambda_S(F_j)}\\
           &=& D(\mu) + \mu(U_j)\frac{\lambda_S(K)}{\lambda_S(F_j)}
  \end{eqnarray}

  \noindent The first equality is Lemma 4.2, the second comes from the
definition of $\mu'_j$ and the last equality uses Lemma 4.3.

  If $\mu$ is optimally dense then by definition $D(\mu) \ge D(\mu_j)$.
So,
$\mu(U_j) = 0$ for all $j$. Hence
  $\mu(\bigcup_j U_j) = 0$. Since $B_j \subset F_j$ for all $j$, $U =
\bigcup_j U_j$ is the set of all packings that are not completely
saturated.
\end{proof}

\section{Tiles with small optimal density}

In this section we will prove Theorem 2.2. We assume from now on that
$S$
is $n$-dimensional hyperbolic space.
    We will need a lemma that follows from much more general ergodic
theory
results of Nevo and Stein. It is an immediate corollary of Theorem 3 in
[NeS].

\begin{lem} If $G$ acts continuously on a compact metric space $X$ such
that there is a Borel probability measure $\mu$ on $X$ that is invariant
and ergodic under this action, then for every function $f \in
L^p(X,\mu)$
($1 < p < \infty$) the following holds for $\mu$-a.e. $x\in X$.

\begin{equation}
\int_X f d\mu = \lim_{n \to \infty} \frac{1}{\lambda_G(\tilde B_n)}
\int_{\tilde B_n} f(gx) \, d\lambda_G(g).
\end{equation}
\end{lem}

\begin{cor}
Let $\mu$ be an ergodic measure in $M_I$. Then the following holds for
$\mu$-a.e. packing $\P \in C_K$.

\begin{equation}
D(\mu) = \lim_{r \to \infty} \frac{\lambda_S(c(\P) \cap
B_r)}{\lambda_S(B_r)}.
\end{equation}
\end{cor}
\begin{proof} By the above lemma, (since $\chi_Z \in L^p(C_K, \mu)$ for
all $p$) for $\mu$-almost every packing $\P$ in $C_K$,
\begin{eqnarray}
D(\mu) &=& \int_{C_K} \chi_Z \, d\mu\\
        &=& \lim_{r \to \infty} \frac{1}{\lambda_G(\tilde B_r)}
\int_{\tilde B_r} \chi_Z(g\P) \, d\lambda_G(g)\\
        &=& \lim_{r \to \infty} \frac{1}{\lambda_S(B_r)} \int_{\tilde
B_r}
\chi_{g^{-1}Z}(\P) \, d\lambda_G(g)
\end{eqnarray}

\noindent If we replace $\pi^{-1}(F_j)$ with $\pi^{-1}(B_r) = \tilde
B_r$
in the proof of equation (13), we get that
\begin{equation}
\int_{\tilde B_n} \chi_{g^{-1}Z}(\P) \, d\lambda_G(g)=
\lambda_S(c(\P) \cap B_r).
\end{equation}
Equations (39) and (40) prove the corollary.
\end{proof}

\noindent Theorem 2.2 will follow almost immediately from this result,
Theorem 3.1 and the following.

\begin{thm} For every $\epsilon >0$ there exists a body $K$ such that
$K$
admits a tiling of $S$ but the optimal density of $K$ is less than
$\epsilon$.
\end{thm}

\noindent {\bf Remark}: It is not yet known whether $K$ can be chosen to
be convex. It is known, however, that there exist convex polygons $P$
which admit tilings of $\H^2$ but such that $D(P) < 1$. For example, the
5-gon with vertices at $i$, $2i$, $i + 1/2$, $i+1$, and $2i+2$ (in the
upperhalf plane model) is one such.\\
\\
To prove this theorem we will construct a body with some indentations
and
protrusions such that in any tiling by $K$ every protrusion in any copy
of
$K$ must fit into an indentation of another copy of $K$. There will be
more indentations than protrusions implying at least, that it does not
admit a tiling of a closed manifold. To show that any invariant measure
in
$C_K$ has small density, we will show that the probability that the origin 
is
in a protrusion is about equal to the probability that the origin is in
any given indentation. Since there are many more indentations than
protrusions, the probability that the origin is in an indentation but
not
in a protrusion will be high, implying that the density is low.

We will need some notation and a lemma (which is similar to a key part
of
the proof of proposition 2 in [BoR]). Let $K$ be a body whose symmetry
group $\Sigma_K$ is trivial. For any $L \subset S$, let $C_K(L) = \{\P
\in
C_K|$ there exists $g \in \hat \P$ such that $\O \in gL\}$.

\begin{lem} Let $K$ be a body with $\Sigma_K = \{e\}$. Let $\mu \in
M_I(K)$ such that $D(\mu) \ne 0$.  If $L$ is a Borel subset of $S$ such
that for every $\P \in C_K(L)$ there exists a unique $g \in \hat \P$
such
that $\O \in gL$ then
\begin{equation}
\mu(C_K(L)) = \mu(C_K(K)) \frac{\lambda_S(L)}{\lambda_S(K)}.
\end{equation}
\end{lem}
\begin{proof}

Let $L_1 = L$ and $L_2 = K$. For $i=1,2$, and $\P \in C_K(L_i)$ let
$\Psi_i(\P)$ be the unique element of $\hat \P$ such that $\O \in
\Psi_i(\P)L_i$.

Note that if for some $i,j \in \{1,2\}$, $\P \in C_K(L_i)$, $g \in G$
and
$g\Psi_i(\P)$ is in the image of $\Psi_j$ then $\O \in
g\Psi_i(\P) K$. Since $\Psi_i(\P) \in \P$, $g\Psi_i(\P) \in g\P$. Hence
$\Psi_j(g\P) = g\Psi_i(\P)$. So if $h$ is in the image of $\Psi_i$ and
$gh$ is in the image of $\Psi_j$ then $\Psi_j^{-1}(gh) =
g\Psi_i^{-1}(h)$.

We define a measure $\lambda'_G$ on $G$ by the following. If $E$ is a
Borel subset of $G$ and there exists a $g \in G$ with $gE$ contained in
the image of $\Psi_i$ for some $i \in \{1,2\}$, then $\lambda'_G(E) =
\mu(\Psi_i^{-1}(gE))$. By the above statement and the $G$-invariance of
$\mu$, this is well-defined. Note that the image of $\Psi_2$ equals
$\pi^{-1}(L_2)^{-1}$ and thus contains an open set (since $K=L_2$
does). So for any Borel subset $E \subset G$ there exists a countable
partition $\{E_i\}_{i=1}^\infty$ of $E$ such that for each $i>0$, there
exists $g_i \in G$ such that $g_i E_i$ is contained in the image of
$\Psi_2$. Then we define $\lambda'_G(E) = \Sigma_{i=1}^\infty
\lambda'_G(E_i)$.

Clearly, $\lambda'_G$ is left $G$-invariant. By the uniqueness of Haar
measure up to scalars (see [Lan]), $\lambda'_G$ must be a scalar
multiple
of $\lambda_G$. Since $D(\mu) = \mu(C_K(K)) \ne 0 $,

\begin{eqnarray}
\frac{\mu(C_K(L_1))}{\mu(C_K(L_2))} &=&
\frac{\lambda'_G(\Psi_1(C_K(L)))}{\lambda'_G(\Psi_2(C_K(K)))}\\
                                     &=&
\frac{\lambda_G(\Psi_1(C_K(L)))}{\lambda_G(\Psi_2(C_K(K)))}\\
                                     &=&
\frac{\lambda_G(\pi^{-1}(L)^{-1})}{\lambda_G(\pi^{-1}(K)^{-1})}\\
                                     &=&
\frac{\lambda_G(\pi^{-1}(L))}{\lambda_G(\pi^{-1}(K))}\\
                                     &=&
\frac{\lambda_S(L)}{\lambda_S(K)}.
\end{eqnarray}

The fourth equation above holds because $\lambda_G$ in inversion
invariant
(i.e. $G$ is unimodular) and the last holds since for any Borel set $E 
\subset S$, $\lambda_S(E) =
\lambda_G(\pi^{-1}(E))$.

\end{proof}

\begin{lem}
If $K, \mu,$ and $L$ satisfy the hypotheses of the above lemma and $L$
is
the disjoint union of $L'$ and $L''$, then both $L'$ and $L''$ satisfy
the
hypotheses of above lemma and $C_K(L)$ is the disjoint union of
$C_K(L')$
and $C_K(L'')$.
\end{lem}
\begin{proof}
This is an easy exercise in understanding the definitions.
\end{proof}

\begin{proof} {\it (of Theorem 5.3)} For simplicity we will only prove the
dimension 2 case. We will identify $S$ with the upperhalf plane model of
$\H^2$ and thus regard it as a subset of the complex plane ${\mathbb C}$
(see [Rat]).

Let $\epsilon>0$ be given. We pick a large positive integer $m$ so that
$\epsilon > 1 - \frac{m-1}{m}$. Let $\delta$ and $\delta'$ be small
positive real numbers with $\delta' < \delta/2 < 1/6$.

Let $R$ be the Euclidean rectangle with vertices $i$, $mi$, $m + mi$,
and
$m+ i$. Let $Q_0$ be the Euclidean rectangle with vertices $\delta + i(1
+
\delta)$, $\delta + i(m - \delta)$, $1 - \delta + i(1 + \delta)$, and $1
-
\delta + i(m - \delta)$. Finally, let $Q'_0$ be the Euclidean rectangle
with vertices $\frac{1}{2} -\delta' + i$, $\frac{1}{2} + \delta' + i$,
$\frac{1}{2} - \delta' + i(1 + \delta)$ and $\frac{1}{2} + \delta' + i(1
+
\delta)$.

For $j = 0,..,m - 1$ let $Q_j = Q_0 + j$ and $Q'_j = Q'_0 + j$. Let $P$
(for protrusion) equal $mQ_0 = \{mz | \, z \in Q\}$ and let $P' = m
Q'_0$. Let $K = R \cup P \cup P' - \cup_j Q_j - \cup Q'_j$.

\begin{figure}[htb]
\begin{center}
\ \psfig{file=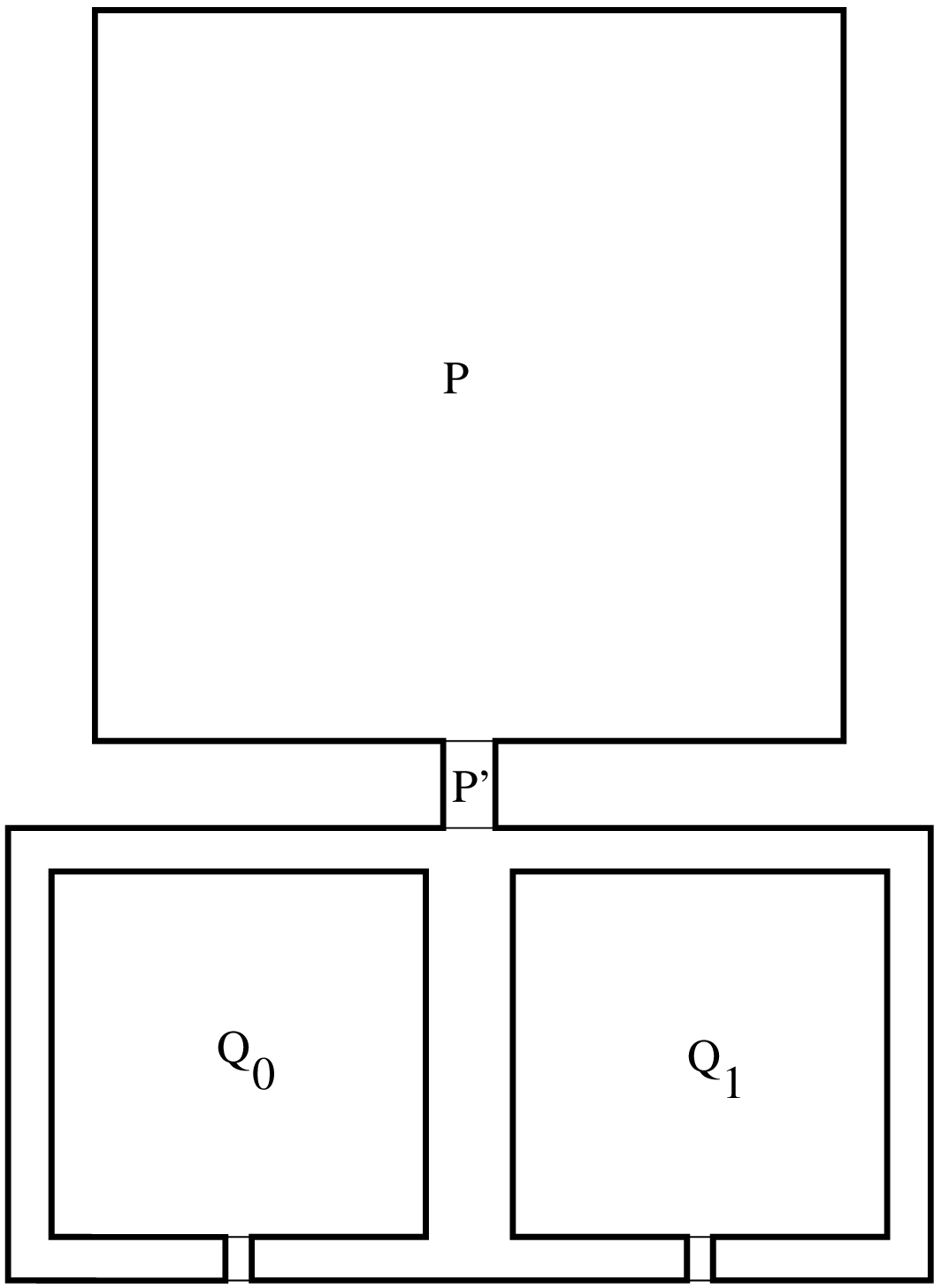,height=2.5in,width=1.87in}
\caption{K}
\label{f:ex}
\end{center}
\end{figure}

In the figure above, $m=2$, $K$ is in heavy outline and $P$, $P'$, $Q_0$, 
$Q_1$, $Q'_0$, and $Q'_1$ are in light outline although $Q'_0$ and $Q'_1$ 
are not labeled.
We let $R' = R \cup P' - \cup_j Q'_j$. Since $ 1 - \frac{m-1}{m} <
\epsilon$ and $\frac{\lambda_S(Q_0)}{\lambda_S(R')} \to \frac{1}{m} $ as
$\delta \to 0$, we can (and do) choose $\delta > 0$ so that
\begin{equation}
1 - (m-1)\left(\frac{\lambda_S(Q_0)}{\lambda_S(R')}\right) < \epsilon.
\end{equation}

We let $\tau_m (z) = z+m$ and $s_m(z) = mz$ for every $z$ in the
upperhalf
plane. We will assume that $\delta' > 0$ has been chosen small enough so
that for any packing $\P \in C_K$ such that the $K$ is in $\P$ and for
any
$gK \in \P$ such that $g \ne s^{-1}_m$, $gK \cap Q_0 = \emptyset$. In
other words, unless the protrusion $gP$ of $gK$ fits into $Q_0$, no part
of $gK$ overlaps $Q_0$. This assumption implies $R'$ is such that for
every $\P \in C_K(R')$, there exists a unique $g \in \hat \P$ such that
$\O \in gR'$. Thus Lemma 5.3 applies.

We now obtain a bound for the optimal density of $K$. Let $\mu \in
M_I(K)$
such that $D(\mu) \ne 0$. By definition of density $D(\mu) =
\mu(C_K(K))$. By Lemma 5.5,

\begin{equation}
D(\mu) = \mu(C_K(K-P)) + \mu(C_K(P)).
\end{equation}

\noindent By Lemma 5.4, $\mu(C_K(P)) = \mu(C_K(Q_j))$ for any $j$. Thus

\begin{equation}
D(\mu) = \mu(C_K(K-P)) + \mu(C_K(Q_0)).
\end{equation}

\noindent By Lemma 5.5,

\begin{eqnarray}
\mu(C_K(R')) &=& \mu(C_K(K-P)) + \Sigma_{j=0}^{m-1} \mu(C_K(Q_j))\\
              &=& \mu(C_K(K-P)) + m\mu(C_K(Q_0)).
\end{eqnarray}

\noindent Equations (49) and (51) imply that

\begin{equation}
D(\mu) = \mu(C_K(R')) - (m-1)\mu(C_K(Q_0)).
\end{equation}

\noindent By Lemma 5.4,
\begin{equation}
\mu(C_K(Q_0)) = \mu(C_K(R'))\frac{\lambda_S(Q_0)}{\lambda_S(R')}.
\end{equation}

\noindent This implies that

\begin{equation}
D(\mu) = \mu(C_K(R'))\left(1 -
(m-1)\frac{\lambda_S(Q_0)}{\lambda_S(R')}\right).
\end{equation}

\noindent Since $\mu(C_K(R')) \le 1$,
\begin{equation}
D(\mu) \le 1 - (m-1)\left(\frac{\lambda_S(Q_0)}{\lambda_S(R')}\right) <
\epsilon.
\end{equation}

\noindent Since $\mu \in M_I$ is arbitrary, $D(K) < \epsilon$. It can
easily be checked that $\{s^i_m \tau^j_m K| \, i, j \in {\mathbb Z}\}$
is
a tiling by $K$.

\end{proof}

\noindent {\bf Remark}: It follows from the above that the space of
tilings by copies of $K$ does not admit a $G$-invariant measure.
However,
this fact can be proven more directly by noting that there is an
equivariant continuous map from the space of tilings by $K$ onto the
space
at infinity of the hyperbolic plane. The latter admits no $G$-invariant
Borel probability measure.

\begin{proof}{\it (of Theorem 2.2)}
Let $\epsilon > 0$. Let $K$ be as in Theorem 5.3. Let $\mu$ be an
optimally dense measure for $K$. Since $\mu$ is ergodic, Corollary
5.2 shows that the set $X_1$ of all $\P \in C_K$ satisfying

\begin{equation}
\lim_{n \to \infty} \frac{\lambda_S(c(\P) \cap B_n)}{\lambda_S(B_n)} =
D(\mu)
\end{equation}

\noindent has $\mu$-measure 1. Let $G'$ be a countable dense subset of
$G$. Then $X_2 = \cap_{g \in G'} \, gX_1$ has $\mu$-measure 1 also. Note 
that for any $\P \in X_2$, $g^{-1}\P \in X_1$ for all $g \in G'$. Hence

\begin{equation}
\lim_{n \to \infty} \frac{\lambda_S(c(\P) \cap B_n(g\O))}{\lambda_S(B_n)} =
D(\mu) < \epsilon
\end{equation}

\noindent for all $g \in G'$. Using the fact that
$G'$ is dense in $G$, it can be shown that the above equation holds for all 
$g \in G$ and for all $\P \in X_2$.

On the other hand, Theorem 3.1 shows that there exists a subset $X_3
\subset C_K$ of full $\mu$-measure such that every packing in $X_3$ is
completely saturated. Then $X_4 = X_2 \cap X_3$ has full $\mu$-measure,
in
particular it is non-empty. Any $\P \in X_4$ satisfies the conclusion of
Theorem 2.2.
\end{proof}


\begin{thebibliography}{xxx}

\bibitem[BoR]{BoR} L. Bowen and C. Radin, \textit{Densest Packings of
Equal Spheres
  in Hyperbolic Space}, preprint, available from
  http://www.ma.utexas.edu/users/radin

\bibitem[FeK]{FeK} G. Fejes Toth and W. Kuperberg, \textit{Packing and
Covering with Convex Sets}, chapter 3.3 in Handbook of Convex Geometry,
ed. P. Gruber and J. Wills, North Holland, Amsterdam, 1993.

\bibitem[FKK]{FKK} G. Fejes Toth, G. Kuperberg, W. Kuperberg,
\textit{Highly
  saturated packings and reduced coverings}, Monatsh. Math. \textbf{125}
(1998),
  no. 2, 127--145.

\bibitem[Gro]{Gro} H. Groemer, \textit{Existensatze fur Lagerungen im
Euclidischen Raum}, Math. Z. \textbf{81} (1963), 260-278.

\bibitem[Lan]{Lan} S. Lang, \textit{Real Analysis}, Addison-Wesley,
Philippines, 1969.

\bibitem[Kup]{Kup} G. Kuperberg, \textit{Notions of denseness},
Geom. Topol. \textbf{4} (2000), 277--292 (electronic),
arXiv:math.MG/9908003.

\bibitem[Nev]{Nev} A. Nevo , \textit{Pointwise ergodic theorems for
radial
averages on simple Lie groups I}, Duke Math. J. \textbf{76} (1994),
113-140.

\bibitem[NeS]{NeS} A. Nevo and E. Stein, \textit{Analogs of Weiner's
ergodic theorems for semisimple groups I}, Annals of Math. \textbf{145}
(1997), 565-595.

\bibitem[Rat]{Rat} J. Ratcliffe, \textit{Foundations of Hyperbolic
Manifolds}, Springer-Verlag, New York, 1994.

\bibitem[RaW]{RaW} C. Radin and M. Wolff, \textit{Space tilings and
local
     isomorphism}, Geometriae Dedicata. \textbf{42} (1992), 355-360.

\end{thebibliography}
\end{document}